\documentclass{cimart}

\usepackage{dsfont}
\usepackage[all]{xy}

\newcommand{\F}{\mathds{F}}
\newcommand{\isoto}{\xrightarrow{\sim}}
\newcommand{\Q}{\mathds{Q}}
\newcommand{\sbul}{{\scriptscriptstyle\bullet}}
\newcommand{\X}{\mathsf{X}}
\newcommand{\Z}{\mathds{Z}}

\DeclareMathOperator{\Gal}{Gal}
\DeclareMathOperator{\inv}{inv}
\DeclareMathOperator{\Maps}{Maps}
\DeclareMathOperator{\Nm}{Nm}

\title[Cup product of inhomogeneous Tate cochains]{Cup product of inhomogeneous Tate cochains, and Galois cohomology of tori over local fields that split over cyclic extensions}

\authors{Mikhail Borovoi}

\authorinfo{Tel Aviv University, Israel}{borovoi@tauex.tau.ac.il}

\abstract{%
    In this note we give formulas for cup product in Tate cohomology in terms of inhomogeneous cochains. Using one of these formulas, for a torus $T$ defined over a non-archimedean local field $K$ and splitting over a cyclic extension of $K$, we compute explicit cocycles representing all cohomology classes in $H^1(K,T)$.
    }



\acknowledgments{%
    The author is grateful to Will Sawin for his answer \cite{Sawin} to the author's MathOverflow question \cite{Borovoi-MO}, and to Tasho Kaletha and James S. Milne for helpful email correspondence. We thank the anonymous referees for their useful comments, which helped us to improve the exposition.
    }

\VOLUME{34}
\ISSUE{1}
\NUMBER{12}
\DOI{https://doi.org/10.46298/cm.17416}
\editinfo{January 28, 2026}{June 23, 2026}{Lenny Fukshansky}

\begin{document}

\section{Introduction}\label{s:intro}

\subsection{Tori over local fields}
\label{ss:tori}
Let $K$ be a non-archimedean local field,
that is, a finite extension of a field of $p$-adic numbers $\Q_p$ for some prime $p$,
or the field of formal power series $\F_q((t))$ over a finite field $\F_q$ of cardinality $q=p^l$
for some natural number $l$.
Let $T$ be a $K$-torus.
For applications in arithmetic geometry, one needs a formula for the Galois cohomology group
\[H^1(K,T):= H^1\big(\Gal(K^s/K), T(K^s)\big)\]
where $K^s$ is a separable closure of $K$.
See Serre \cite[Section I.2.2]{Serre} for the definition of group cohomology.

Let $L/K$ be a finite Galois extension splitting  $T$.
Consider $G=\Gal(L/K)$, its Galois group,
and the cocharacter group of our torus
\[X={\sf X}_*(T):={\rm Hom}({\mathds G}_{m,L},T_L).\]
Then $G$ naturally acts on $X$.
We have
\[H^1(K,T)=H^1(L/K,T):= H^1(G,T(L));\]
see Sansuc \cite[(1.9.2), page 19]{Sansuc}.
We have a canonical isomorphism
\[T(L)=X\otimes L^\times\]
where tensor product is taken over $\Z$.
Therefore, we wish to compute
\[H^1(G, X\otimes L^\times).\]

According to Tate \cite[Theorem on page 717]{Tate-Nagoya},
there is a canonical isomorphism
\begin{equation}\label{e:tau}
 H^{-1}(G,X)\isoto  H^1(G,X\otimes L^\times).
\end{equation}
Recall that, denoting  by $X_G$ the group of coinvariants of $G$ in $X$:
\[X_G=X/\langle x-g\cdot x\ |\ g\in G, x\in X\rangle,\]
the group $H^{-1}(G,X)$ is the subgroup of $X_G$
consisting of the classes of elements $x\in X$ with $\sum_{g\in G}\, g\cdot x=0$.
Let
\[u_{L/K}\in H^2(G,L^\times)\]
denote the {\em fundamental class}; see, for instance, Harari \cite[Definition 9.1]{Harari}.
The isomorphism \eqref{e:tau}
is given by cup product with $u_{L/K}$.

Tate's formula \eqref{e:tau} describes the finite abelian group
$H^1(K,T)=H^1(L/K,T)$
up to a canonical isomorphism.
However, for applications, in particular, for twisting,
it would be helpful to have not only the isomorphism class of this abelian group,
but also an explicit formula for the isomorphism \eqref{e:tau}.
For this end, we need formulas for the fundamental class and for cup product.

In the general case of finite Galois extension $L/K$ of local fields,
the fundamental class is  mysterious.
In the case where the extension $L/K$ is unramified,
and hence cyclic (that is, $G$ is cyclic),
there is a formula for the fundamental class $u_{L/K}$
in Milne's notes \cite{Milne-CFT}.

In the answer \cite{Sawin} to a MathOverflow question of the author, Will Sawin
guessed and partly proved a formula for $u_{L/K}$
in the case when $L/K$ is cyclic (possibly ramified).
We prove Sawin's formula in full using Milne's answer to the author's question \cite{Borovoi-MO}:

\begin{theorem}
\label{t:Sawin}
Let  $L/K$ be a \textbf{cyclic} extension of degree $n$ of non-archimedean local fields.
Choose a generator $\sigma$ of the cyclic group (of order $n$) $G=\Gal(L/K)$.
Let $e_\sigma\in K^\times$ be an element such that
\[ (e_\sigma,L/K)=\sigma\]
where
\[ e\mapsto (e,L/K)\colon\ K^\times \to K^\times/\Nm (L^\times) \isoto G\]
is the local reciprocity homomorphism; see \cite[Section VI.2.2]{CF}.
Then the fundamental class $u_{L/K}\in H^2(G,L^\times)$
is represented by the following  cocycle $a\colon G\times G\to L^\times$:
\[\text{for \ $0\le i,j<n$,}\qquad a(\sigma^i,\sigma^j)=
\begin{cases}
1 &\text{if \ \ } i+j< n,\\
e_\sigma &\text{if \ \ }i+j\ge n.
\end{cases}
\]
\end{theorem}

\subsection{The homogeneous and inhomogeneous complexes}
\label{ss:cup}

Our formulas for $H^{-1}(G,X)$, for $u_{L/K}\in  H^2(G,L^\times)$,  and for $H^1(G, X\otimes L^\times)$,
are written in terms of {\em inhomogeneous} cocycles.
On the other hand, the formulas for cup product in \cite[Section IV.7]{CF}
are written in terms of the {\em homogeneous} cocycles.
Below we write what the corresponding homogeneous and inhomogeneous complexes are.

Let $G$ be a finite group and $A$ be a $G$-module.
We consider the {\em homogeneous} cochain complex
$F^\sbul$ for $A$ (see \cite[Appendix E]{BK}):
{\small
\begin{equation}
\label{e:c-homogeneous}
\begin{aligned}
\xymatrix@C=3mm@R=3mm{
\dots F^{-2}( G,A)\ar[r]^-{d^{-1}}  &F^{-1}( G,A)\ar[r]^-{d^{0}}  &F^{0}( G,A)\ar[r]^-{d^{1}}  &F^{1}( G,A)\ar[r]^-{d^2}&F^{2}( G,A)\dots \\
 \Maps_ G( G^2,A)\ar@{=}[u] &\Maps_ G( G,A)\ar@{=}[u] &\Maps_ G( G,A)\ar@{=}[u] &\Maps_ G( G^2,A)\ar@{=}[u]&\Maps_ G( G^3,A)\ar@{=}[u]
}
\end{aligned}
\end{equation}
}

We also consider the {\em inhomogeneous} cochain complex $C^\sbul$ for $A$ (see \cite[Appendix E]{BK}):
{\small
\begin{equation}
\label{e:c-inhomogeneous}
\begin{aligned}
\xymatrix@C=6.8mm@R=3mm{
\dots C^{-2}( G,A)\ar[r]^-{d^{-1}}  &C^{-1}( G,A)\ar[r]^-{d^0}  &C^{0}( G,A)\ar[r]^-{d^1}  &C^{1}( G,A)\ar[r]^-{d^2}  &C^{2}( G,A)\dots\\
\Maps( G,A)\ar@{=}[u] &A\ar@{=}[u] &A\ar@{=}[u] &\Maps( G,A)\ar@{=}[u] &\Maps( G^2,A)\ar@{=}[u].
}
\end{aligned}
\end{equation}
}

These two complexes $F^\sbul$ and $C^\sbul$ are canonically isomorphic; see Section \ref{s:cup} below.
For $p\in\Z$, the elements of $F^p(G,A)$ are called  homogeneous $p$-cochains,
and the elements of $C^p(G,A)$ are called inhomogeneous $p$-cochains.

Let $A$ and $B$ be two $G$-modules.
For $p,q\in \Z$, there is a canonical {\em cup product pairing}
\begin{equation}\label{e:cup-intro}
\cup\colon F^p(G,A)\times F^q(G,B)\to F^{p+q}(G,A\otimes B);
\end{equation}
see \cite[Section IV.7]{CF} or Section \ref{s:cup} below
(here $p$ is  unrelated to the residue characteristic of $K$).
It induces a cup product pairing
\begin{equation}\label{e:cup-intro-inhom}
\cup\colon C^p(G,A)\times C^q(G,B)\to C^{p+q}(G,A\otimes B).
\end{equation}

When $p,q\ge0$, the ``homogeneous'' formula for the  cup product pairing \eqref{e:cup-intro}
is simple, and so is the  ``inhomogeneous'' formula for \eqref{e:cup-intro-inhom}.
For instance, if
\[
    c_1' \in C^1(G,A)=\Maps(G, A)
    \quad \text{and} \quad
    c_1'' \in C^1(G,B)=\Maps(G, B)
\]
are two inhomogeneous cochains, then
for their cup product,
\[
    c_2=c_1'\cup c_1''\in C^2(G,A)=\Maps(G\times G,A\otimes B),
\]
we have
\[ c_2(h_1,h_2)=c_1'(h_1)\otimes h_1\!\cdot\! c''_2(h_2)\quad\ \text{for}\ \, h_1,h_2\in G.\]
However, when $p$ and $q$ have different signs, the corresponding homogeneous formulas contain summation,
and so do the inhomogeneous formulas.
In Section \ref{s:cup}, using the canonical isomorphism
between the homogeneous and the  inhomogeneous complexes,
we give formulas for the cup product pairing in terms of inhomogeneous cochains \eqref{e:cup-intro-inhom}
for all $p$ and $q$.

\subsection{Tate's isomorphism made explicit.}
Using Theorem \ref{t:Sawin} and one of our formulas for cup product in terms of inhomogeneous cochains,
we explicitly compute the isomorphism \eqref{e:tau} in the cyclic case:

\begin{theorem}\label{t:-1}
Let $T$  be a  torus over a non-archimedean local field $K$, and assume that
$T$ splits over a \textbf{cyclic} extension $L$ of $K$ of degree $n$.
Let $G=\Gal(L/K)$, and let $\sigma\in G$ and $e_\sigma\in K^\times$ be
as in Theorem \ref{t:Sawin}.
Let $x\in X=\X_*(T)$ be a cocharacter such that $[x]\in X_G$
is contained in $H^{-1}(G,X)$.
Then the image of $[x]\in H^{-1}(G,X)$ in $H^1(G,X\otimes L^\times)$
under the isomorphism \eqref{e:tau} is the class of the cocycle
\[z_x\colon \,G\to X\,\otimes\,L^\times, \quad\ \sigma^g\,\longmapsto\
    \sum_{t=1}^{g}\sigma^t\cdot x\,\otimes\,e_\sigma\quad  \text{for}\ \,0\le g<n.\]
\end{theorem}

The plan of the rest of the note is as follows.
In Section \ref{s:cup} we write the formulas for cup product
in terms of homogeneous and inhomogeneous cocycles.
In Section \ref{s:fundamental} we prove Theorem \ref{t:Sawin}.
In Section \ref{s:torus} we prove Theorem \ref{t:-1}.
In Appendix \ref{app:cup} we provide details of calculations in Section \ref{s:cup}.

\section{Cup product in terms of inhomogeneous cocycles}
\label{s:cup}

In this section, $G$ is an arbitrary finite group and $A$ is a $G$-module.

We consider the homogeneous cochain complex
\eqref{e:c-homogeneous} for $A$, which we denote by
$F^\sbul$.
We also consider the inhomogeneous complex
\eqref{e:c-inhomogeneous} for $A$, which we denote by $C^\sbul$.
For details, see \cite[Appendix E]{BK}.

These two complexes $F^\sbul$ and $C^\sbul$ are canonically isomorphic,
with the isomorphism given as follows.
For a homogeneous cochain $f\in F^n( G,A)$
with corresponding inhomogeneous cochain $c\in C^n( G,A)$, we have
\begin{align}
c(h_1,\dots,h_n)&=f(  1,  h_1,  h_1 h_2, \dots,  h_1\cdots h_n ),\quad n\ge 0
\label{e:c+}\\
f(g_0,\dots,g_n)&=g_0c(  g_0^{-1}g_1,  g_1^{-1}g_2,  g_2^{-1}g_3, \dots,  g_{n-1}^{-1}g_n ),\quad\ n\ge0
\label{e:f+}\\
c(h_1,\dots,h_{-n-1})&=f(  1,h_1,  h_1 h_2, \dots,  h_1\cdots h_{-n-1} ),\quad n<0
\label{e:c-}\\
f(g_1,\dots,g_{-n})&=g_1c(  g_1^{-1}g_2,  g_2^{-1}g_3, \dots,  g_{-n-1}^{-1}g_{-n} ),\quad\ n<0.
\label{e:f-}
\end{align}

We give formulas for cup product of cochains, homogeneous and inhomogeneous.
The formulas for homogeneous cochains follow immediately from
\cite[Section IV.7, formulas for $\varphi_{p,q}$ on page 107]{CF}.
We deduce formulas for inhomogeneous cochains from these formulas
using formulas \eqref{e:c+}--\eqref{e:f-}.
Details of calculations are given in Appendix \ref{app:cup}.

Let $A$ and $B$ be two $G$-modules.
For $m,n\in\Z$, consider homogeneous cochains
\[f'_m\in F^m(G,A),\quad f''_n\in F^n(G,B),\quad f_{m+n}=f'_m\cup f''_n\in F^{m+n}(G,A\otimes B).\]
Consider the corresponding inhomogeneous cochains
\[c'_m\in C^m(G,A),\quad c''_n\in C^n(G,B).\]
We find formulas for  the inhomogeneous cochain $c_{m+n}\in C^{m+n}(G,A\otimes B)$
corresponding to $f_{m+n}$.
We write then
\[ c_{m+n}= c'_m\cup c''_n .\]

Let  $p,q\ge 0$. We take $m=p$, $n=q$. Then we have
\begin{align}
\label{e:f++}
f_{p+q}&(g_0,\dots,g_{p+q})=
   f'_p(g_0,\dots,g_p)\,\otimes\,f''_q(g_p,\dots,g_{p+q})
\\
\label{e:c++}
c_{p+q}&(h_1,\dots,h_{p+q})\\
&=c'_p(h_1,\dots,h_p)\,\otimes\, (h_1h_2\cdots h_p)\cdot c''_q(h_{p+1},  \dots,  h_{p+q}).
\notag
\end{align}

Let  $p,q\ge1$. We take $m=-p$, $n=-q$. Then we have
\begin{align}
\label{e:f--}
f_{-p-q}&(g_0,\dots,g_{p+q-1})=f'_{-p}(g_0,\dots,g_{p-1})\,\otimes\,f''_{-q}(g_{p},\dots,g_{p+q-1}),\\
\label{e:c--}
c_{-p-q}&(h_1,\dots,h_{p+q-1})\\
&=c'_{-p}(h_1,\dots,h_{p-1}) \otimes  (h_1h_2\cdots h_p)\cdot c''_{-q}(h_{p+1},  \dots,  h_{p+q-1}).
\notag
\end{align}

Let $p\ge 0$, $q\ge 1$. We take $m=p$, $n=-p-q$. Then we have
\begin{align}
\label{e:f+--}
&f_{-q}(\,g_1,\ \dots,\ g_q\,)=\sum_s f'_{p}(g_1,\,s_1,\ \dots,\ s_p\,)
    \,\otimes\,f''_{-p-q}(s_{p},\ \dots,\ s_1, g_1,\ \dots,\ g_q\,)\\
\label{e:c+--}
&c_{-q}(h_1,\ \dots,\ h_{q-1})\\
&\quad=\sum_t c'_p(\,t_1,\ \dots,\ t_p\,)\,\otimes\,(t_1\cdots t_p)\cdot
c''_{-p-q}(\,t_p^{-1},\ \dots,\ t_1^{-1}, h_1,\ \dots,\ h_{q-1}\,).
\notag
\end{align}

Let $p\ge 0$, $q\ge 1$. We take $m=-p-q$, $n=p$. Then we have
\begin{align}
\label{e:f--+}
&f_{-q}(\,g_1,\ \dots,\ g_q\,)=\sum_s f'_{-p-q}(g_1,\ \dots,\ g_q, s_1, \ \dots,\ s_p)
\,\otimes\,f''_p(s_p,\ \dots,\ s_1, g_q).\\
\label{e:c--+}
&c_{-q}(h_1,\ \dots,\ h_{q-1})\\
&\ \ =\sum_t c'_{-p-q}(h_1,\ \dots,\ h_{q-1},t_1,\ \dots,\ t_p)
\,\otimes\, (h_1h_2\cdots h_{q-1}t_1\cdots t_p)\cdot c''_p(t_p^{-1},\ \dots,\ t_1^{-1}).
\notag
\end{align}

Let $p\ge 0$, $q\ge 1$. We take $m=p+q$, $n=-q$. Then we have
\begin{align}
\label{e:f++-}
&f_p(g_0,\dots,g_p)=\sum_s f'_{p+q}(g_0,\dots,g_p,s_1,\dots, s_{q})\,\otimes\,f''_{-q}(s_{q},\dots, s_1 )\\
\label{e:c++-}
&c_p(h_1,\dots,h_p)\\
&\quad=\sum_t  c'_{p+q}\big(  h_1, \dots,h_p,t_1,\dots,t_q)
   \,\otimes\,(h_1\cdots h_p t_1\cdots t_q) \cdot c''_{-q}(t_q^{-1},\dots, t_2^{-1})
\notag
\end{align}

Let $p\ge 0$, $q\ge 1$.
We take $m=-q$, $n=p+q$. Then we have
\begin{align}
\label{e:f-++}
&f_p(g_0,\ \dots,\ g_p)=\sum_s f'_{-q}(s_1,\ \dots,\ s_q)\,\otimes\,f''_{p+q}(s_q,\ \dots,\ s_1, g_0,\ \dots,\ g_p)\\
\label{e:c-++}
&c_p( h_1, \ \dots,\ h_p)\\
&\quad=\sum_t t_1\cdot c'_{-q}(t_2,\ \dots,\ t_q)
   \,\otimes\,(t_1\cdots t_q)\cdot c''_{p+q}(t_q^{-1},\ \dots,\ t_1^{-1},\ h_1,\ \dots,\ h_p).
\notag
\end{align}
\noindent In these formulas, the sums are taken over $s_i$ (resp. $t_i$) running independently through $G$.

\section{The fundamental class of a cyclic extension}
\label{s:fundamental}
In this section we prove Theorem \ref{t:Sawin}.
Let $L/K$ be a finite Galois extension of non-archimedean local fields with {\em cyclic} Galois group
$ G=\Gal(L/K)$ of order $n$.
Choose a generator $\sigma\in G$.
Let $e_\sigma\in K^\times$ be as in Subsection \ref{ss:tori}.
We write $[e_\sigma]$ for the image of $e_\sigma$
in $H^{-1}(G,L^\times):= K^\times/\Nm(L^\times)$.

Consider the  homomorphism
$$\chi\colon  G\to \frac1n\Z/\Z\quad\ \text{such that}\ \,\chi(\sigma)=\frac1n.$$
Then $\chi\in H^1( G,\Q/\Z)$.
Consider
$$\delta\chi\in H^2( G,\Z),$$
the image of $\chi$ under the connecting map
$$\delta\colon H^1( G,\frac1n\Z/\Z)\to H^2( G,\Z).$$
Then $\delta\chi$ is represented by the 2-cocycle
$$ b\in Z^2( G,\Z),\quad\ b\colon G\times G\to \Z$$
defined as follows (compare Milne \cite[page 102, above Proposition III.1.9]{Milne-CFT}).
We represent an element of $ G\times G$ as $(\sigma^i,\sigma^j)$
with $i,j\in\Z$, $0\le i,j<n$.
Then
$$
b(\sigma^i,\sigma^j)=
\begin{cases}
0 &\text{if}\quad i+j<n\\
1 &\text{if} \quad i+j\ge n.
\end{cases}
$$

For $a$ as in Theorem \ref{t:Sawin}, we have
$$ a=  b\cup e_\sigma $$
where
$e_\sigma\in Z^0( G, L^\times)=K^\times$ \,
and $\cup$ denotes cup product of cocycles.
It follows that $a$ is a 2-cocycle, $a\in Z^2( G, L^\times)$.
Write $[a]\in H^2( G,L^\times)$ for the class of $a$.
Write
$$[e_\sigma]=e_\sigma\cdot \Nm(L^\times)\in K^\times/\Nm(L^\times)=\hat H^0( G, L^\times).$$
Then
$$[a]= [b]\cup [e_\sigma]=\delta\chi\cup [e_\sigma]$$
where $\cup$ denotes cup product of Tate cohomology classes.
Consider the canonical isomorphism
$${\rm inv}_{L/K}\colon H^2( G,L^\times)\isoto \frac1n\Z/\Z.$$
We have
$${\rm inv}_{L/K}( \delta\chi\cup [e_\sigma])=\chi(\sigma);$$
see Harari \cite[Proposition 9.3]{Harari}.
By construction, we have $\chi(\sigma)=\frac1n$.
Thus,
\[\inv_{L/K}[a]=\inv_{L/K}(\delta\chi\cup[e_\sigma])=\chi(\sigma)=\frac1n.\]
By the definition of the fundamental class $u_{L/K}$, it follows that $[a]=u_{L/K}$,
which completes the proof of Theorem \ref{t:Sawin}.

\section{Explicit cocycles for a torus}
\label{s:torus}

In this section we prove Theorem \ref{t:-1}.
We sometimes write $z$ for $z_x$.
We compute
\[ z=x\cup a\in Z^{-1}(G,X)\otimes  Z^2(G,L^\times) \]
where
\begin{align*}
&z\in Z^1(G, X\otimes L^\times)\subseteq\Maps(G,X\otimes L^\times),\\
&a\in Z^2(G,L^\times)\subseteq\Maps(G\times G,L^\times),\\
&x\in Z^{-1}(G,X)\subseteq X.
\end{align*}
Here $[a]=u_{L/K}$, the fundamental class, and the 2-cocycle $a$ is as in Theorem \ref{t:Sawin}
 and its proof, namely, $a=b\cup e_\sigma$.

Formula \eqref{e:c-++} for $p=1$, $q=1$ reads:
\[ z(g)=\sum_{t\in G} t\cdot x\,\otimes\, t\cdot a(t^{-1}, g)\]
for  $g\in G$. Since $a$ takes the values $e_\sigma\in K^\times$ and $1$ only, we have
$t\cdot a(t^{-1}, g) = a(t^{-1}, g)$.
It remains to compute the sum
\[ \sum_{t\in G} t\cdot x\,\otimes\, a(t^{-1}, g).\]

Recall that $\sigma$ is a generator of the cyclic group $G$.
By abuse of notation, we write $\sigma^g$ for $g$ and $\sigma^t$ for $t$.
Now $g,t\in \Z/n\Z$.

We wish to compute
\[
z(\sigma^g)=\sum_{t\in\Z/n\Z} \sigma^t\cdot x\,\otimes\, a(\sigma^{-t},\sigma^g).
\]
For this end, we compute
\[a(\sigma^{-t},\sigma^g)= a(\sigma^{n-t},\sigma^g)\quad\ \text{where}\ \, 0\le g<n,\quad  0<t\le n. \]

If $0<t\le g<n$, then
\[ 0< n-t <n, \quad 0\le g<n,\quad (n-t)+g= n+(g-t)\ge n \ \ \text{(because $t\le g$)},\]
whence $a(\sigma^{n-t}, \sigma^g)=e_\sigma$.

If $t>g\ge 0,\ \, t\le n$, then
\[0\le n-t <n, \quad 0\le g <n,\quad (n-t)+g=n+(g-t)<n \ \ \text{(because $t> g$)},\]
whence $a(\sigma^{n-t}, \sigma^g)=1$.

We conclude that
\begin{equation}\label{e:g<n}
z_x(\sigma^g)=\sum_{t=1}^{g}\sigma^t\!\cdot\! x\,\otimes\,e_\sigma
\end{equation}
for $ 0\le g<n$.
This completes the proof of  Theorem \ref{t:-1}.

\begin{remark}
Recall that
\[Z^{-1}(G,X)=\big\{  x\in X\ \,\big|\ \,\sum_{t=1}^n \sigma^t\!\cdot\! x=0 \big\}.\]
It follows that for every $x\in Z^{-1}(G,X)$ and every $k\in \Z$ we have
\[\sum_{t=k+1}^{k+n}\sigma^t\!\cdot\!x= \sigma^k\cdot \Big(\sum_{s=1}^{n}\sigma^s\!\cdot\!x\Big)=\sigma^k\cdot 0=0.\]
Thus, equality \eqref{e:g<n} holds for all $g\ge 0$.
Now we can check explicitly that $z_x$ in \eqref{e:g<n} is indeed a 1-cocycle.
For every $g_1,g_2\ge0$ we have
\begin{align*}
z_x(\sigma^{g_1+g_2})&= \sum_{t=1}^{g_1+g_2}\sigma^t\!\cdot\! x\,\otimes\,e_\sigma\\
&=\sum_{t=1}^{g_1}\sigma^t\!\cdot\! x\,\otimes\,e_\sigma
+ \sum_{t=g_1+1}^{g_1+g_2}\sigma^t\!\cdot\! x\,\otimes\,e_\sigma\\
&=\sum_{t=1}^{g_1}\sigma^t\!\cdot\! x\,\otimes\,e_\sigma
+\sigma^{g_1}\!\cdot\!\Big(\sum_{s=1}^{g_2}\sigma^s\!\cdot\! x\,\otimes\,e_\sigma\Big)\\
&=z_x(\sigma^{g_1})+ \sigma^{g_1}\!\cdot\! z_x(\sigma^{g_2}),
\end{align*}
as desired.
\end{remark}

\appendix

\section{Cup product: calculations}
\label{app:cup}

\subsection{}
We derive formula \eqref{e:c++} from \eqref{e:f++}.
For $p,q\ge 0$, set $m=p$, $n=q$.
Then by \cite[Section~IV.7, page 107]{CF}
we have \eqref{e:f++}:
\begin{align*}
f_{p+q}&(g_0,\dots,g_{p+q})=
   f'_p(g_0,\dots,g_p)\,\otimes\,f''_q(g_p,\dots,g_{p+q}).
\end{align*}
Substituting \eqref{e:f+} in the right-hand side,
we obtain
\[f_{p+q}(g_0,\dots,g_{p+q})=g_0\cdot c'_p(  g_0^{-1}g_1, \dots,  g_{p-1}^{-1}  g_p )
\,\otimes\,g_p\cdot c''_q(  g_p^{-1}g_{p+1}, \dots,  g_{p+q-1}^{-1}  g_{p+q} ).
\]
Substituting
\[g_0=1,\ \, g_1=h_1,\ \, g_2=h_1 h_2,\ \dots,
   \ g_p=h_1 h_2\cdots h_p,\ \dots,\ g_{p+q}=h_1 h_2\cdots h_{p+q} ,\]
we obtain
\begin{align*}
&g_0^{-1}g_1=h_1,\ \, g_1^{-1}g_2=h_2, \ \dots,\ g_{p-1}^{-1}  g_{p}=h_{p}, \\
&g_{p}^{-1} g_{p+1}=h_{p+1},\ \dots,\  g_{p+q-1}^{-1}  g_{p+q}=h_{p+q},
\end{align*}
whence
\begin{align*}
f_{p+q}(1,h_1,h_1h_2,\ \dots,\ h_1h_2\cdots h_{p+q})
=c'_p(h_1,\dots,h_p) \otimes  (h_1h_2\cdots h_p)\cdot c''_q(h_{p+1},  \dots,  h_{p+q}),
\end{align*}
and applying \eqref{e:c+}  to the left-hand side, we obtain \eqref{e:c++}:
\begin{align*}
c_{p+q}(h_1,\dots,h_{p+q})
=c'_p(h_1,\dots,h_p) \otimes  (h_1h_2\cdots h_p)\cdot c''_q(h_{p+1},  \dots,  h_{p+q}).
\end{align*}

\subsection{}
We derive formula \eqref{e:c--} from \eqref{e:f--}.
For $p,q\ge1$, set $m=-p$, $n=-q$. Then by \cite[Section IV.7, page 107]{CF} we have \eqref{e:f--}:
\begin{align*}
f_{-p-q}&(g_0,\dots,g_{p+q-1})=f'_{-p}(g_0,\dots,g_{p-1})\,\otimes\,f''_{-q}(g_{p},\dots,g_{p+q-1}).
\end{align*}
Substituting \eqref{e:f-} in the right-hand side,
we obtain
\begin{multline*}
f_{-p-q}(g_0,\dots,g_{p+q-1})\\
=g_0\cdot  c'_{-p}(g_0^{-1}g_1,g_1^{-1}g_2,\dots,g_{p-2}^{-1}g_{p-1})
\,\otimes\,g_p\cdot c''_{-q}(g_p^{-1}g_{p+1}, \dots, g_{p+q-2}^{-1}g_{p+q-1}).
\end{multline*}
Substituting
\[g_0=1,\ \, g_1=h_1,\ \, g_2=h_1 h_2,\ \dots,
   \ g_p=h_1 h_2\cdots h_p,\ \dots,\ g_{p+q-1}=h_1 h_2\cdots h_{p+q-1} ,\]
we obtain
\begin{align*}
&g_0^{-1}g_1=h_1,\ \, g_1^{-1}g_2=h_2,  \ \dots,\ g_{p-2}^{-1}  g_{p-1}=h_{p-1}, \\
&g_{p}^{-1} g_{p+1}=h_{p+1},\ \dots,\   g_{p+q-2}^{-1}  g_{p+q-1}=h_{p+q-1},
\end{align*}
whence
\begin{multline*}
f_{-p-q}(1,h_1,h_1h_2,\dots,h_1h_2\cdots h_{p+q-1})\\
=c'_{-p}(h_1,\dots,h_{p-1}) \otimes  (h_1h_2\cdots h_p)
    \cdot c''_{-q}(h_{p+1},  \dots,  h_{p+q-1}),
\end{multline*}
and applying \eqref{e:c-} to the left-hand side, we obtain \eqref{e:c--}:
\begin{align*}
c_{-p-q}(h_1,\dots,h_{p+q-1})
=c'_{-p}(h_1,\dots,h_{p-1}) \otimes  (h_1h_2\cdots h_p)\cdot c''_{-q}(h_{p+1},  \dots,  h_{p+q-1}).
\end{align*}

\subsection{}
From now on we assume that  $p\ge 0$, $q\ge 1$.
We derive formula \eqref{e:c+--} from \eqref{e:f+--}.

We take $m=p$, $n=-p-q$. Then by \cite[Section IV.7, page 107]{CF} we have \eqref{e:f+--}:
\begin{align*}
f_{-q}(\,g_1,\ \dots,\ g_q\,)=\sum_s f'_{p}(g_1,\,s_1,\ \dots,\ s_p\,)
    \,\otimes\,f''_{-p-q}(s_{p},\ \dots,\ s_1, g_1,\ \dots,\ g_q\,)
\end{align*}
where  $s_i$  run independently through $G$.
Substituting \eqref{e:f+} and  \eqref{e:f-} in the right-hand side,
we obtain
\begin{multline*}
f_{-q}(\,g_1,\ \dots,\ g_q\,)=\sum_s\Big(
g_1\cdot c'_p(g_1^{-1}s_1, s_1^{-1}s_2,\ \dots,\ s_{p-1}^{-1}s_p)\\
\,\otimes\, s_p\cdot c''_{-p-q}(s_p^{-1}s_{p-1},\, \dots,\, s_2^{-1}s_1, s_1^{-1}g_1, g_1^{-1}g_2,
\, \dots,\, g_{q-1}^{-1}g_q\Big).
\end{multline*}
Substituting
\begin{align*}
&s_1=t_1,\, s_2=t_1t_2,\ \dots,\ s_p=t_1\cdots t_p,\\
&g_1=1,\ g_2=h_1,\, g_3=h_1 h_2, \ \dots,\ g_q=h_1\cdots h_{q-1},
\end{align*}
we obtain
\begin{align*}
&g_1^{-1}s_1=t_1,\ s_1^{-1} s_2=t_2, \ \dots,\ s_{p-1}^{-1}s_p=t_p,\\
&s_p^{-1}s_{p-1}=t_p^{-1},\ \dots,\  s_2^{-1} s_1=t_2^{-1},\ s_1^{-1} g_1=t_1^{-1},\\
&g_1^{-1}g_2=h_1, \ g_2^{-1} g_3=h_2,\ \dots,\ g_{q-1}^{-1} g_q=h_{q-1},
\end{align*}
whence
\begin{multline*}
f_{-q}(\,1, h_1, h_1h_2,\ \dots,\ h_1h_2\cdots h_{q-1}\,)\\
=\sum_t c'_p(\,t_1,\ \dots,\ t_p\,)\,\otimes\,(t_1\cdots t_p)\cdot
c''_{-p-q}(\,t_p^{-1},\ \dots,\ t_1^{-1}, h_1,\ \dots,\ h_{q-1}\,),
\end{multline*}
and applying \eqref{e:c-} to the left-hand side, we obtain \eqref{e:c+--}:
\begin{align*}
c_{-q}(h_1,\ \dots,\ h_{q-1})
=\sum_t c'_p(\,t_1,\ \dots,\ t_p\,)\,\otimes\,(t_1\cdots t_p)\cdot
c''_{-p-q}(\,t_p^{-1},\ \dots,\ t_1^{-1}, h_1,\ \dots,\ h_{q-1}\,)
\end{align*}
where  $t_i$  run independently through $G$.

\subsection{}
We derive formula \eqref{e:c--+} from \eqref{e:f--+}.
We take $m=-p-q$, $n=p$. Then by \cite[Section~IV.7, page 107]{CF}  we have \eqref{e:f--+}:
\begin{align*}
f_{-q}(\,g_1,\ \dots,\ g_q\,)=\sum_s f'_{-p-q}(g_1,\ \dots,\ g_q, s_1, \ \dots,\ s_p)
\,\otimes\,f''_p(s_p,\ \dots,\ s_1, g_q)
\end{align*}
where  $s_i$  run independently through $G$.
Substituting \eqref{e:f-} and \eqref{e:f+} in the right-hand side,
we obtain
\begin{multline*}
f_{-q}(\,g_1,\ \dots,\ g_q\,)\\
=\sum_s g_1\cdot c'_{-p-q}(g_1^{-1}g_2,\ \dots,\ g_{q-2}^{-1} g_{q-1},
\ g_q^{-1} s_1,\ s_1^{-1} s_2, \ \dots, s_{p-1}^{-1}s_p) \\
\otimes\,
s_p\cdot c''_p(s_p^{-1}s_{p-1},\ \dots,\ s_2^{-1}s_1, s_1^{-1}g_q).
\end{multline*}
Substituting
\begin{gather*}
g_1=1,\ g_2=h_1,\ g_3=h_1h_2,\ \dots,\ g_q=h_1h_2\cdots h_{q-1},\\
s_1=h_1h_2\cdots h_{q-1}t_1,\ \dots,\ s_p=h_1h_2\cdots h_{q-1}t_1\cdots t_p ,
\end{gather*}
we obtain
\begin{align*}
&g_1^{-1}g_2=h_1,\ \dots,\ g_{q-1}^{-1} g_q=h_{q-1},\\
&g_q^{-1}s_1=t_1,\  s_1^{-1}s_2=t_2,\ \dots,\ s_{p-1}^{-1}s_p=t_p,\\
&s_p^{-1}s_{p-1}=t_p^{-1},\ \dots,\  s_2^{-1}s_1=t_2^{-1},\ s_1^{-1}g_q=t_1^{-1},
\end{align*}
whence
\begin{multline*}
f_{-q}(\,1,\ h_1,\ h_1h_2,\ \dots,\ h_1h_2\cdots h_{q-1}\,)\\
=\sum_t c'_{-p-q}(h_1,\ \dots,\ h_{q-1},t_1,\ \dots,\ t_p)
\,\otimes\,(h_1h_2\cdots h_{q-1}t_1\cdots t_p)\cdot c''_p(t_p^{-1},\ \dots,\ t_1^{-1}),
\end{multline*}
and applying \eqref{e:c-} to the left-hand side, we obtain \eqref{e:c--+}:
\begin{multline*}
c_{-q}(\,h_1,\ \dots,\ h_{q-1}\,)\\
\sum_t c'_{-p-q}(h_1,\ \dots,\ h_{q-1},t_1,\ \dots,\ t_p)
\,\otimes\,(h_1h_2\cdots h_{q-1}t_1\cdots t_p)\cdot c''_p(t_p^{-1},\ \dots,\ t_1^{-1})
\end{multline*}
where  $t_i$  run independently through $G$.

\subsection{}
We derive formula \eqref{e:c++-} from \eqref{e:f++-}.
We take $m=p+q$, $n=-q$. Then by \cite[Section IV.7, page 107]{CF}
we have \eqref{e:f++-}:
\begin{align*}
f_p(g_0,\dots,g_p)=\sum_s f'_{p+q}(g_0,\dots,g_p,s_1,\dots, s_{q})\,\otimes\,f''_{-q}(s_{q},\dots, s_1 )
\end{align*}
where  $s_i$  run independently through $G$.
Substituting \eqref{e:f+} and  \eqref{e:f-} in the right-hand side,
we obtain
\begin{multline*}
f_p(\,g_0,\ \dots,\ g_p\,)\\
=\sum_s\Big( g_0\cdot c'_{p+q}\big(  g_0^{-1}g_1,\ \dots,\
    g_{p-1}^{-1} g_p,\ g_p^{-1}s_1,\ s_1^{-1}s_2, \dots,\ s_{q-1}^{-1} s_{q} \big)\\
   \,\otimes\,s_q \cdot c''_{-q}(s_q^{-1}s_{q-1},\ \dots,\  s_2^{-1}s_1)\Big).
\end{multline*}
Substituting
\begin{align*}
&g_0=1,\ \, g_1=h_1,\ \, g_2=h_1 h_2,\ \dots,
   \ g_p=h_1 h_2\cdots h_p,\\
&s_1=h_1h_2\cdots h_p t_1,\ s_2=h_1h_2\cdots h_p t_1 t_2,\ \dots,\ s_q=h_1h_2\cdots h_p t_1\cdots t_q ,
\end{align*}
we obtain
\begin{align*}
&g_0^{-1}g_1=h_1,\ \, g_1^{-1}g_2=h_2,
  \ \dots,\ \,g_{p-1}^{-1} g_p=h_p,\\
&g_p^{-1} s_1=t_1,\ \, s_1^{-1}s_2=t_2,\ \dots,\   s_{q-1}^{-1}s_{q}=t_{q},\\
&s_q^{-1} s_{q-1}=t_q^{-1},\ \dots,\ s_2^{-1} s_1=t_2^{-1},
\end{align*}
whence
\begin{multline*}
f_p(1,h_1,h_1h_2,\ \dots,\ h_1h_2\cdots h_p)\\
\sum_t  c'_{p+q}\big(\,  h_1,\  \dots,\ h_p,t_1,\ \dots,\ t_q\,\big)
   \,\otimes\,(h_1\cdots h_p t_1\cdots t_q) \cdot c''_{-q}(t_q^{-1},\ \dots,\  t_2^{-1}),
\end{multline*}
and applying \eqref{e:c+} to the left-hand side, we obtain \eqref{e:c++-}:
\begin{align*}
c_p&(h_1,\dots,h_p) \\
&=\sum_t  c'_{p+q}\big(\,  h_1,\  \dots,\ h_p,t_1,\ \dots,\ t_q\,\big)
   \,\otimes\,(h_1\cdots h_p t_1\cdots t_q) \cdot c''_{-q}(t_q^{-1},\ \dots,\  t_2^{-1})
\end{align*}
where  $t_i$  run independently through $G$.

\subsection{}
We derive formula \eqref{e:c-++} from \eqref{e:f-++}.
We take $m=-q$, $n=p+q$. Then by \cite[Section IV.7, page 107]{CF}
we have \eqref{e:f-++}:
\begin{align*}
f_p(g_0,\ \dots,\ g_p)=\sum_s f'_{-q}(s_1,\ \dots,\ s_q)\,\otimes\,f''_{p+q}(s_q,\ \dots,\ s_1, g_0,\ \dots,\ g_p)
\end{align*}
where  $s_i$  run independently through $G$.
Substituting \eqref{e:f-} and  \eqref{e:f+} in the right-hand side,
we obtain
\begin{multline*}
f_p(g_0,\ \dots,\ g_p)\\
=\sum_s s_1\cdot c'_{-q}(s_1^{-1}s_2,\ \dots,\ s_{q-1}^{-1} s_q)
   \,\otimes\,s_q\cdot c''_{p+q}(s_q^{-1} s_{q-1},\ \dots,\ s_1^{-1} g_0,\ g_0^{-1}g_1,\ \dots,\ g_{p-1}^{-1} g_p).
\end{multline*}
Substituting
\begin{align*}
&g_0=1,\ \, g_1=h_1,\ \, g_2=h_1 h_2,\ \dots,  \ g_p=h_1 h_2\cdots h_p,\\
&s_1=t_1,\ \, s_2=t_1t_2,\ \dots,\ s_q=t_1\cdots t_q ,
\end{align*}
we obtain
\begin{align*}
&s_1^{-1}s_2=t_2,\ \dots,\ s_{q-1}^{-1} s_q=t_q,\\
&s_q^{-1}s_{q-1}=t_q^{-1},\ \dots,\ s_1^{-1}g_0=t_1^{-1} ,\\
&g_0^{-1}g_1=h_1, \ \dots,\ g_{p-1}^{-1}g_p=h_p,
\end{align*}
whence
\begin{multline*}
f_p(1,\ h_1, \ \dots,\ h_1\cdots h_p)\\
=\sum_t t_1\cdot c'_{-q}(t_2,\ \dots,\ t_q)
   \,\otimes\,(t_1\cdots t_q)\cdot c''_{p+q}(t_q^{-1},\ \dots,\ t_1^{-1},\ h_1,\ \dots,\ h_p),
\end{multline*}
and applying \eqref{e:c+} to the left-hand side, we obtain \eqref{e:c-++}:
\begin{align*}
c_p( h_1, \ \dots,\ h_p)=
\sum_t t_1\cdot c'_{-q}(t_2,\ \dots,\ t_q)
   \,\otimes\,(t_1\cdots t_q)\cdot c''_{p+q}(t_q^{-1},\ \dots,\ t_1^{-1},\ h_1,\ \dots,\ h_p)
\end{align*}
where  $t_i$  run independently through $G$.


\end{document}